\documentclass[12pt]{amsart}
\usepackage{amsmath, amsthm, amsfonts, amssymb}

\newtheorem{theorem}{Theorem}[section]
\newtheorem{lemma}[theorem]{Lemma}

\theoremstyle{definition}

\newtheorem{example}[theorem]{Example}

\theoremstyle{remark}

\numberwithin{equation}{section}

\newcommand{\mC}{\ensuremath{\mathbb{C}}}

\newcommand{\mN}{\ensuremath{\mathbb{N}}}

 \allowdisplaybreaks

\begin{document}
\title{Some Variations on an Extension of Carath$\acute{E}$odory's Normality Criterion}

\author[K. S. Charak]{Kuldeep Singh Charak}
\address{
\begin{tabular}{lll}
&Kuldeep Singh Charak\\
&Department of Mathematics\\
&University of Jammu\\
&Jammu-180 006\\ 
&India\\
\end{tabular}}
\email{kscharak7@rediffmail.com}

\author[R. Kumar]{Rahul Kumar}
\address{
\begin{tabular}{lll}
&Rahul Kumar\\
&Department of Mathematics\\
&University of Jammu\\
&Jammu-180 006\\
&India
\end{tabular}}
\email{rktkp5@gmail.com}
 
\begin{abstract}
 In this note we give some variations on an extension of Carath$\acute{e}$odory's normality criterion due to Grahl and Nevo.
\end{abstract}

\renewcommand{\thefootnote}{\fnsymbol{footnote}}
\footnotetext{2010 {\it Mathematics Subject Classification} 30D45, 30D30.}
\footnotetext{{\it Keywords and phrases}. Normal families, Carath$\acute{e}$odory's normality criterion, Spherical metric.}
\footnotetext{The work of the first author is partially supported by Mathematical Research Impact Centric Support (MATRICS) grant, File No. MTR/2018/000446, by the Science and Engineering Research Board (SERB), Department of Science and Technology (DST), Government of India.\\
The work of the second author is supported by the CSIR India.}

\maketitle

\section{Introduction and Main Results}

\medskip

By $\mathcal{M}(D)$ we shall denote the class of all meromorphic functions on a domain $D$ in $\mC _{\infty},$ where  $\mC _{\infty}$ has two metrics $\chi$ and $\sigma$ known as the chordal and spherical metrics.  $\chi$ and $\sigma$ are equivalent metrics and hence interchangeable. A subfamily $\mathcal{F}$ of $\mathcal{M}(D)$ is said to be normal on $D$ if  we can extract a subsequence from each given sequence in $ \mathcal{F}$ that converges locally uniformly on $D$ to a limit function which is either  meromorphic on $D$ or identically equal to $ \infty.$  The famous Montel's Fundamental Normality Test (see \cite{Schiff}) states that {\it if each member of $\mathcal{F}$ omits three distinct fixed values in  $\mC_\infty,$ then $\mathcal{F}$ is normal in $D.$} 

\smallskip

Carath$\acute{e}$odory\cite{Car} proved the following extension of Montel's Fundamental Normality Test:\\

\textbf{Theorem A:}
{\it Let $\mathcal{F}$ be a subfamily of $\mathcal{M}(\mathbb{D})$ and suppose that there exists an $\epsilon>0$ such that each $f\in\mathcal{F}$ omits three distinct values $a_f,b_f,c_f \in \mC_\infty$ satisfying
$$min\{\sigma(a_f,b_f),\sigma(b_f,c_f),\sigma(c_f,a_f)\}\geq\epsilon.$$
Then $\mathcal{F}$ is normal in $\mathbb{D}.$}

\smallskip

Grahl and Nevo\cite{Grahl} generalized Theorem A as:\\

\textbf{Theorem B:}
{\it Let $\mathcal{F}$ be a subfamily of $\mathcal{M}(\mathbb{D})$ and $\epsilon>0.$ Suppose that for each $f\in\mathcal{F}$ there exist three functions $a_f,b_f,c_f\in\mathcal{M}(\mathbb{D})\cup\{\infty\}$ such that $f$ omits the functions $a_f,b_f,c_f$ in $\mathbb{D}$ and
$$min\{\sigma(a_f(z),b_f(z)),\sigma(b_f(z),c_f(z)),\sigma(c_f(z),a_f(z))\}\geq\epsilon$$
for all $z\in \mathbb{D}.$ Then $\mathcal{F}$ is normal in $\mathbb{D}.$}

\smallskip

The purpose of this paper is to present some variations on Theorem B. Precisely, we consider the following variants:
\begin{itemize}
\item[(i)] $f-a_{jf}$ and $g-a_{jg}$ share $0$ for each pair of functions $f, g\in \mathcal{F}$;
\item[(ii)] $f^{(k)}-a_{jf}$ and $g^{(k)}-a_{jg}$ share $0$ for each pair of functions $f, g\in \mathcal{F}$;
\item[(iii)] For each $f\in \mathcal{F}$, the zeros of $f-a_{1f}$ are of multiplicity at least four and those of $f-a_{jf}$ for $j=2,3,$ are of multiplicity at least three.
\end{itemize}
\begin{theorem}\label{thm1}
Let $\mathcal{F}\subseteq\mathcal{M}(D)$ be such that for each $f\in\mathcal{F}$ there exist three functions $a_{1f},a_{2f},a_{3f}\in\mathcal{M}(D)\cup\{\infty\}$ satisfying
\begin{equation}
\label{eq1}
\sigma(a_{jf}(z), \ a_{kf}(z))\geq \epsilon ~\mbox{for}~ j\neq k ~\mbox{and}~ z\in D, \mbox{ for some } \epsilon >0.
\end{equation}
If for each pair of functions $f ~\mbox{and}~ g$ in $\mathcal{F}, \  f-a_{jf} \mbox{ and } g-a_{jg}$ share $0$ in $D, ~(j= 1,2,3).$ Then $\mathcal{F}$ is normal in $D.$
\end{theorem}

\begin{example} \label{ex1}
Consider the family $\mathcal{F}= \{f_n:n\in\mathbb{N}\},$ where $$f_n(z)= nz$$ on the open unit disk $\mathbb{D}$ and let $a_{1f_n}(z)= nz, ~a_{2f_n}(z)= nz^2, ~a_{3f_n}(z)= nz^3.$ Clearly, for every $f,g\in\mathcal{F}, f-a_{jf} ~\mbox{and}~ g-a_{jg} ~(j= 1,2,3)$ share the value $0$ in $\mathbb{D}.$ However, the family $\mathcal{F}$ is not normal in $\mathbb{D}.$ Note that the functions $a_{jf}, j=1, 2, 3$ do not satisfy $(\ref{eq1}).$
\end{example}
 
\begin{example} \label{ex2}
Consider the family $\mathcal{F}= \{f_n:n\in\mathbb{N}\},$ where $f_n(z)= \tan{nz}$ on the open unit disk $\mathbb{D}$ and let $a_{1f_n}(z)= i, ~a_{2f_n}(z)= -i.$
Take $\epsilon= \sigma(i,-i).$ Clearly, for every $f,g\in\mathcal{F}, f-a_{jf} ~\mbox{and}~ g-a_{jg} ~(j= 1,2)$ share the value $0$ in $\mathbb{D}.$ Also $\sigma(a_{1f_n}, a_{2f_n})\geq\epsilon.$ However, the family $\mathcal{F}$ is not normal in $\mathbb{D}.$  Thus the number of functions $a_{jf}$ in Theorem \ref{thm1} cannot be reduced to two.
\end{example}

\begin{theorem}\label{thm1.1}
Let $\mathcal{F}\subseteq\mathcal{M}(D)$ be such that all zeros of each $f\in \mathcal{F}$ are of  multiplicity at least $k+1$ and that for each $f \in \mathcal{F}$ there exist three functions $a_{1f},a_{2f},a_{3f}\in\mathcal{M}(D)\cup\{\infty\}$ satisfying 
\begin{equation}
\label{eq2}
\sigma(a_{jf}(z), \ a_{lf}(z))\geq \epsilon ~\mbox{for}~ j\neq l ~\mbox{and}~ z\in D, \mbox{ for some } \epsilon >0.
\end{equation}
If for each pair of functions $f ~\mbox{and}~ g$ in $\mathcal{F}, \  f^{(k)}-a_{jf} \mbox{ and } g^{(k)}-a_{jg}$ share $0$ in $D, ~(j= 1,2,3).$ Then $\mathcal{F}$ is normal in $D.$
\end{theorem}
The condition, ``zeros of each $f \in \mathcal{F}$ have multiplicity at least $k+1$" in Theorem \ref{thm1.1} cannot be dropped:
\begin{example}
Consider the family $\mathcal{F}= \{f_n:n\in\mathbb{N}\},$ where $$f_n(z)= nz^k$$ on the open unit disk $\mathbb{D}.$ Let 
$$a_{1f_n}(z)= 1/2, ~a_{2f_n}(z)= 1/3, ~a_{3f_n}(z)= 1/4$$ 
and take 
$$\epsilon= \min\left\{\sigma(1/2, 1/3), \sigma(1/2, 1/4), \sigma(1/3, 1/4)\right\}.$$
Clearly, for every $f,g\in\mathcal{F}, f^{(k)}-a_{jf} ~\mbox{and}~ g^{(k)}-a_{jg} ~(j= 1,2,3)$ share the value $0$ in $\mathbb{D}.$ Also $\sigma(a_{jf_n}, a_{kf_n})\geq\epsilon.$ However, the family $\mathcal{F}$ is not normal in $\mathbb{D}.$ 
\end{example}

\begin{theorem}\label{thm2} Let $\mathcal{F}\subseteq\mathcal{M}(D)$ be such that for each $f\in\mathcal{F}$ there exist functions $a_{1f},a_{2f},a_{3f}\in\mathcal{M}(D)\cup\{\infty\}$ such that
\begin{itemize}
\item[(i)] $\sigma(a_{jf}(z),a_{kf}(z))\geq \epsilon ~\mbox{for}~ j\neq k ~\mbox{and}~ z\in D \mbox{ for some } \epsilon >0;$
\item[(ii)] each zero of $f-a_{1f}$ has multiplicity at least four;
\item[(iii)] each zero of $f-a_{jf}, (j=2,3)$ has multiplicity at least three.
\end{itemize}
 Then $\mathcal{F}$ is normal in $D.$
\end{theorem}
The following example shows that the condition, ``$\sigma(a_{jf}(z),a_{kf}(z))\geq \epsilon ~\mbox{for}~ j\neq k ~\mbox{and}~ z\in D$"  in Theorem \ref{thm2} cannot be dropped:

\begin{example}
 Consider $\mathcal{F}=\{f_n\}, ~\mbox{where}~ f_n(z)= nz^4+nz^3$ on $\mathbb{D}.$\\
Let $a_{1f_n}(z)= nz^3, a_{2f_n}(z)= nz^4 ~\mbox{and}~ a_{3f_n}(z)= n(z-1/2)^3+nz^4+nz^3.$ Then each zero of $f_n-a_{1f_n}$ has multiplicity at least four and each zero of $f_n-a_{jf_n}, (j=2,3)$ has multiplicity at least three but $\mathcal{F}$ is not normal on $\mathbb{D}.$
\end{example}

\begin{theorem}\label{thm3}
Let $\mathcal{F}$ be a subfamily of $\mathcal{M}(D)$ and let $h\in\mathcal{M}(D).$ Then $\mathcal{F}$ is normal in $D$ if and only if each $z_0\in D$ has a neighbourhood $D(z_0; r)$ such that each $f\in\mathcal{F}$ satisfies:
$$\sigma(f(z),h(z))>\epsilon_1  ~\forall z\in D(z_0;r) ~\mbox{or}~ \sigma(f(z),h(z))<\epsilon_2  ~\forall z\in D(z_0;r)$$
 for some  $ \epsilon_1, \epsilon_2 \mbox{ with } 0<\epsilon_1 <\epsilon_2 <\pi.$
\end{theorem}

\section{Preparation for the proofs of theorems }
Let $f$ be meromorphic function in $\mC$ and $ a\in \mC_\infty.$ Then $a$ is called totally ramified value of $f$  if $f-a$ has no simple zeros.
Following result known as {\it Nevanlinna's Theorem } (see \cite{Bergweiler}) plays a crucial role in the proof of Theorem \ref{thm2}:

\smallskip

\textbf{Nevanlinna's Theorem:}
{\it Let $f$ be a non-constant meromorphic function,  $a_1,...,a_q\in \mC_\infty \mbox{ and } m_1,...,m_q \in \mN.$ Suppose that all $a_j$-points of $f$ have multiplicity  at least $m_j,  \mbox{ for }~ j=1,...,q.$ Then 
$$\sum\limits_{j=1}^{q}(1-\frac{1}{m_j}) \leq 2.$$
}
If $f$ does not assume  the value $a_j$ at all, then  we take  $m_j=\infty.$ 

\medskip

For $p\in\mathbb{N}$ and $j= 1,\ldots,p$ we define the projections $\pi_j:(\mathcal{M}(D))^p\to\mathcal{M}(D)$ as:
$$\pi_j(f_1,\ldots,f_p):=f_j ~\mbox{for}~  (f_1,\ldots,f_p)\in(\mathcal{M}(D))^p.$$

\begin{lemma}(\cite{Grahl})\label{lemma01} Let $p\in\mathbb{N}$ and $\mathcal{F}\subseteq(\mathcal{M}(\mathbb{D}))^p.$ Assume that there exists $j_0\in\{1,2,\ldots,p\}$ such that the family $\pi_{j_0}(\mathcal{F})$ is not normal at $z_0\in \mathbb{D}.$ Then there exist sequences $(f_n)_n= ((f_{1,n},\ldots,f_{p,n}))_n\subseteq\mathcal{F}, (z_n)_n\subseteq\mathbb{D} ~\mbox{and}~ (\rho_n)_n\subseteq (0, ~1)$ such that $\lim\limits_{n\to\infty}z_n= z_0, \lim\limits_{n\to\infty}\rho_n= 0$ and such that for all $j= 1,2,\ldots,p$ the sequences $(g_{j,n})_n$ defined by $$g_{j,n}(\zeta):= f_{j,n}(z_n+\rho_n\zeta)$$
converges to functions $g_j\in\mathcal{M}(\mathbb{C})\cup\{\infty\}$ locally uniformly in $\mathbb{C}$ with respect to the spherical metric, where at least one of the functions $g_1,g_2,\ldots,g_p$ is non-constant.
\end{lemma}

\begin{lemma}(\cite{Grahl})\label{lemma1} Let $\mathcal{G}\subseteq(\mathcal{M}(\mathbb{D}))^2$ be a family of pairs of meromorphic functions in $\mathbb{D}$ and $\epsilon > 0.$ Assume that 
$$\sigma(a(z),b(z))\geq\epsilon,~\mbox{for all}~ (a,b)\in\mathcal{G} ~\mbox{and all}~ z\in\mathbb{D}.$$
Then the families $\{a:(a,b)\in\mathcal{G}\}$ and $\{b:(a,b)\in\mathcal{G}\}$ are normal in $\mathbb{D}.$
\end{lemma}                      

\begin{lemma}(\cite{Chen})\label{lemmaC}
Let $\mathcal{F}\subseteq\mathcal{M}(D),$ all of whose zeros have multiplicity at least $k+1,$  where $k$ is a positive integer. If $\mathcal{F}^{(k)}= \{f^{(k)}:f\in\mathcal{F}\}$ is normal in  $D,$ then  $\mathcal{F}$ is also normal in $D.$ 
\end{lemma}

\begin{lemma}\label{lemma3}
 If $\mathcal{G}\subseteq (\mathcal{M}(D))^2 ~\mbox{and}~  0<\epsilon<\pi.$ If 
$$\sigma(a(z),b(z))<\epsilon ~\mbox{for all}~ (a,b)\in\mathcal{G}$$
and  $\{b:(a,b)\in\mathcal{G}\}$ is normal in $D.$ Then $\{a:(a,b)\in\mathcal{G}\}$ is also normal in $D.$
\end{lemma}
\textbf{Proof of Lemma \ref{lemma3}:} 
Suppose on the contrary that $\{a:(a,b)\in\mathcal{G}\}$ is not normal at $z_0\in D.$ That is, the first projection $\pi_1(\mathcal{G})$ is not normal at $z_0.$ Then by Lemma \ref{lemma01} there exist sequences $\{(a_n,b_n)\}\subseteq\mathcal{G}, \{z_n\}\subseteq D ~\mbox{and}~ \{\rho_n\}\subset (0, 1)$ such that $z_n\to z_0, \rho_n\to 0 ~\mbox{as}~ n\to\infty,$ and the sequences $\{g_n\} ~\mbox{and}~ \{h_n\}$ defined as 
$$g_n(\zeta):= a_n(z_n+\rho_n\zeta) ~\mbox{and}~ h_n(\zeta):= b_n(z_n+\rho_n\zeta)$$
converging locally uniformly in $\mC$ (with respect to the spherical metric) to functions $g ~\mbox{and}~ h$ respectively,  where at least one of $g$ or $h$ is non-constant. Since $\{b:(a,b)\in\mathcal{G}\}$ is normal in $D, h$ is constant on $\mathbb{C}, ~\mbox{say}, h(\zeta)= w_0, \forall\zeta\in\mathbb{C}.$ Since
$$\sigma(a_n(z_n+\rho_n\zeta),b_n(z_n+\rho_n\zeta))<\epsilon,$$
it follows that $\sigma(g(\zeta), w_0)\leq\epsilon, ~\mbox{for all}~ \zeta\in\mathbb{C}$ which is a contradiction since $g$ is non-constant.  $\Box$

\section{Proofs of Theorems}
\textbf{Proof of Theorem \ref{thm1}:}  Let $$\mathcal{G}:= \{(f,a_{1f},a_{2f},a_{3f}):f\in\mathcal{F}\}$$ 
 Suppose $\mathcal{F}$ is not normal at $z_0\in D.$ That is, the first projection $\pi_1(\mathcal{G})$ is not normal at $z_0.$ Then by Lemma \ref{lemma01} there exist sequences $\{(f_n,a_{1f_n},a_{2f_n},a_{3f_n})\}\subseteq\mathcal{G}, \{z_n\}\subseteq D ~\mbox{and}~ \{\rho_n\}\subset (0, 1)$ such that  
$z_n\to z_0, \rho_n\to 0 ~\mbox{as}~ n\to\infty$ and the sequences $\{g_n\}, \{A_{1n}\},$\\
 $\{A_{2n}\}, \{A_{3n}\}$ defined  by
$$g_n(\zeta):= f_n(z_n+\rho_n\zeta),$$
$$A_{1n}(\zeta):= a_{1f_n}(z_n+\rho_n\zeta), A_{2n}(\zeta):= a_{2f_n}(z_n+\rho_n\zeta), A_{3n}(\zeta):= a_{3f_n}(z_n+\rho_n\zeta) $$
converges spherically locally uniformly in $\mathbb{C}$ to functions $g, A, B, C\in\mathcal{M}(\mathbb{C})\cup\{\infty\}$ respectively, all of which are not constant.\\
Since $$\sigma(a_{jf_n}(z),a_{kf_n}(z))\geq\epsilon ~\mbox{for}~ j\neq k ~\mbox{and}~ z\in D ~\mbox{and}~ n\in\mathbb{N}$$
it follows that $\{a_{if_n}\}, ~(i=1,2,3)$ forms a normal family on $D$ and hence the limit functions $A,B,C$ are constant and consequently $g$ is non-constant.\\
{\it Claim:} $g$ omits at least two of the values ${A, B, C.}$\\
Suppose on the contrary $g$ assumes two values $B ~\mbox{and}~ C.$ Let $g(\zeta_0)= B ~\mbox{and}~ g(\zeta_0^*)=C ~\mbox{for some}~ \zeta_0, \zeta_0^*\in \mathbb{C}.$ By Hurwitz's theorem there exist sequences $\zeta_n\to \zeta_0$
and $\zeta_n^*\to \zeta_0^*$ such that 
$$f_n(z_n+\rho_n\zeta_n)-a_{2f_n}(z_n+\rho_n\zeta_n)=0 ~\mbox{and}~ f_n(z_n+\rho_n\zeta_n^*)-a_{3f_n}(z_n+\rho_n\zeta_n^*)=0.$$
By hypothesis, for any positive integer $m,$
$$f_m(z_n+\rho_n\zeta_n)-a_{2f_m}(z_n+\rho_n\zeta_n)=0 ~\mbox{and}~ f_m(z_n+\rho_n\zeta_n^*)-a_{3f_m}(z_n+\rho_n\zeta_n^*)=0.$$
Fix $m$ and take $n\to\infty$, we have 
$$f_m(z_0)-a_{2f_m}(z_0)= 0 ~\mbox{and}~ f_m(z_0)-a_{3f_m}(z_0)=0.$$
which contradicts the fact that $\sigma(a_{2f_m}(z_0), a_{3f_m}(z_0))\geq\epsilon.$

Without loss of generality, we assume $g$ omits $B ~\mbox{and}~ C$ and so $g(\zeta)$ assume $A$ infinitely many times. Let $\zeta^{(1)} ~\mbox{and}~ \zeta^{(2)}$ be two distinct zeros of $g(\zeta)-A.$ Choose  $\delta>0$ small enough such that $D(\zeta^{(1)},\delta)\cap D(\zeta^{(2)},\delta)= \emptyset .$\\
Since $g$ is non-constant, by Hurwitz's theorem there are points $\zeta_n^{(1)}\in D(\zeta^{(1)},\delta)$ and $\zeta_n^{(2)}\in D(\zeta^{(2)},\delta)$ such that for sufficiently large  $n$
$$f_n(z_n+\rho_n\zeta_n^{(1)})-a_{1f_n}(z_n+\rho_n\zeta_n^{(1)})=0 ~\mbox{and}~ f_n(z_n+\rho_n\zeta_n^{(2)})-a_{1f_n}(z_n+\rho_n\zeta_n^{(2)})=0.$$
By hypothesis, for any positive integer $m,$
$$f_m(z_n+\rho_n\zeta_n^{(1)})-a_{1f_m}(z_n+\rho_n\zeta_n^{(1)})=0 ~\mbox{and}~ f_m(z_n+\rho_n\zeta_n^{(2)})-a_{1f_m}(z_n+\rho_n\zeta_n^{(2)})=0.$$
Fix $m,$ take $n\to\infty,$ and in view of the fact that $z_n+\rho_n\zeta_n^{(1)}\to z_0 ~\mbox{and}~ z_n+\rho_n\zeta_n^{(2)}\to z_0,$ we have
$$f_m(z_0)-a_{1f_m}(z_0)= 0.$$
Since the zeros of $f_m-a_{1f_m}$ have no accumulation point, we have 
$$z_n+\rho_n\zeta_n^{(1)}= z_0 ~\mbox{and}~ z_n+\rho_n\zeta_n^{(2)}= z_0$$
and so $\zeta_n^{(1)}= \zeta_n^{(2)},$ which is a contradiction. $\Box$

\smallskip

\textbf{Proof of Theorem \ref{thm1.1}:} By the similar arguments as in the proof of Theorem \ref{thm1}, we conclude that $\mathcal{F}^{(k)}= \{f^{(k)}:f\in\mathcal{F}\}$ is normal on  $D$ and hence by Lemma \ref{lemmaC}, the normality of $\mathcal{F}$ on $D$ follows.  $\Box$

\smallskip

\textbf{Proof of Theorem \ref{thm2}:}  Let $$\mathcal{G}:= \{(f,a_{1f},a_{2f},a_{3f}):f\in\mathcal{F}\}$$ 
 Suppose $\mathcal{F}$ is not normal at $z_0\in D.$ That is, the first projection $\pi_1(\mathcal{G})$ is not normal at $z_0.$ Then by Lemma \ref{lemma01} there exist sequences $\{(f_n,a_{1f_n},a_{2f_n},a_{3f_n})\}\subseteq\mathcal{G}, \{z_n\}\subseteq D ~\mbox{and}~ \{\rho_n\}\subset (0, 1)$ such that  
\begin{itemize}
\item[(a)] $z_n\to z_0, \rho_n\to 0 ~\mbox{as}~ n\to\infty; $
\item[(b)] all zeros of $f_n-a_{1f_n}$ has multiplicity at least four ;
\item[(b)] all zeros of $f_n-a_{jf_n}, (j= 2,3)$ has multiplicity at least three ;
\item[(c)] $\sigma(a_{jf_n}(z),a_{kf_n}(z))\geq\epsilon ~\mbox{for}~ j\neq k ~\mbox{and}~ z\in D ~\mbox{and}~ n\in\mathbb{N};$ and
\item[(d)] the sequences $\{g_n\}, \{A_{1n}\}, \{A_{2n}\}, \{A_{3n}\}$ defined by
$$g_n(\zeta):= f_n(z_n+\rho_n\zeta),$$
$$A_{1n}(\zeta):= a_{1f_n}(z_n+\rho_n\zeta), A_{2n}(\zeta):= a_{2f_n}(z_n+\rho_n\zeta), A_{3n}(\zeta):= a_{3f_n}(z_n+\rho_n\zeta) $$
converges spherically locally uniformly in $\mathbb{C}$ to functions $g, A, B, C\in\mathcal{M}(\mathbb{C})\cup\{\infty\}$ respectively, all of which are not constant.
\end{itemize}
Now by $(c)$ it follows that $\{a_{jf_n}\}, (j= 1,2,3)$ forms a normal family in $D$ and hence the limit functions $A,B,C$ are constant and consequently $g$ is non-constant.\\
Let $h_n(\zeta):= g_n(\zeta)-A_{1n}(\zeta).$ Then $h_n(\zeta)\to g(\zeta)-A.$ Now by $(b)$ and Argument principle, it follows that each zero of $g(\zeta)-A$ has multiplicity at least four.  That is, A-points of g are of multiplicity at least four. Similarly, we find that B-points and C-points of $g$ are of multiplicity at least three. Since $g$ is non-constant, we arrive at a contradiction by Nevanlinna's Theorem.  $\Box$

\smallskip

\textbf{Proof of Theorem \ref{thm3}:} Suppose $\mathcal{F}$ is normal in $D$ and suppose on the contrary that there is a neighborhood $D(z_0,r)$ of $z_0,$ a sequence $\{f_n\} ~\mbox{in}~ \mathcal{F}$ two sequences $\{z_n\} ~\mbox{and}~ \{w_n\} $ of complex numbers with $z_n\to z_0 ~\mbox{and}~ w_n\to z_0$ such that
$$\sigma(f_n(z_n),h(z_n))\leq\epsilon_1 ~\mbox{and}~ \sigma(f_n(w_n),h(w_n))\geq\epsilon_2$$
Since $\{f_n\}$ is normal on $D,$ there is a subsequence (we still denote it by $\{f_n\}$) such that $\{f_n\}$ converges locally uniformly to $f$ with respect to the spherical metric. Thus

\begin{eqnarray*}
\epsilon_1 &\geq& \lim\limits_{n\to\infty} \sigma(f_n(z_n),h(z_n))\\
           &=& \sigma(f(z_0),h(z_0))\\
					 &=&	\lim\limits_{n\to\infty} \sigma(f_n(w_n),h(w_n))\\
					&\geq& \epsilon_2,
\end{eqnarray*}
which is a contradiction. 

\smallskip

Conversely, suppose that each $z_0\in D$ has a neighborhood with the given property. Define $h_n(z):= h(z) ~\forall n.$ Then $\{h_n\}$ is a normal family in $D$ and $h_n \to h$ locally uniformly with respect to the spherical metric. Then for each $n$
$$\sigma(f_n(z),h_n(z))>\epsilon_1 ~\forall z\in D(z_0,r) ~\mbox{or}~ \sigma(f_n(z),h_n(z))<\epsilon_2 ~\forall z\in D(z_0,r). $$
In the first case $\{f_n\}$  is normal in $D(z_0,r)$ by Lemma \ref{lemma1} and in the second case $\{f_n\}$ is normal in $D(z_0,r)$ by Lemma \ref{lemma3}. $\Box$ 

\section{Normality of Families of Rational Maps}

Beardon, Minda and Short(\cite{Beardon}) proved the following theorem:
\begin{theorem}
\label{thmM}
For any family  $\mathcal{F}$ of M$\ddot{o}$bius maps and any domain $D$ in $\mathbb{C}_\infty$, the following are equivalent:
\begin{itemize}
\item[(a)] The family $\mathcal{F}$ is normal in $D.$
\item[(b)] For any set $A= \{u,v,w\}$ of distinct points in $\mathbb{C}_\infty,$ each point $z_0$ in $D$ has a neighbourhood in which each $f$ in $\mathcal{F}$ omits at least two of the values in $A.$
\item[(c)] There exists a set $A= \{u,v,w\}$ of distinct points in $\mathbb{C}_\infty$ such that  each point $z_0$ in $D$ has a neighbourhood in which each $f$ in $\mathcal{F}$ omits at least two of the values in $A.$
\end{itemize} 
\end{theorem}
Whether Theorem \ref{thmM} holds for rational functions of higher degree, is the objective of this section. Let $\mathcal{R}_m$ denote a family of rational maps of degree $d:0\leq d\leq m, \  m\in \mathbb{N}.$ Suppose $\mathcal{R}_m$ is normal in $D.$ Then consider the set $B:=\{u, v, w, x \}$, where $x\in \mathbb{C}$ is different from other elements of the set. Then each $z_0 \in D$ has a neighbourhood in which each $R\in \mathcal{R}_m$ omits at least three of the four values in $B$ showing that $(a)$ implies $(b)$ holds for the family $\mathcal{R}_m$. Also $(b) \Rightarrow (c)$ holds trivially for $\mathcal{R}_m.$ Now we prove that $(c)\Rightarrow (a)$ also holds for $\mathcal{R}_m:$\\
Suppose $(c)$ holds. Let $\{R_p\}$ be any sequence in $\mathcal{R}_m$ and $z_0\in D.$ Then $z_0$ has a neighborhood $U_0$ in which each $R_p$  omits at least two of the  values in $A.$ We consider the subsets  $\{R_{p_j}\}, \  \{R_{p_k}\}$ and $\{R_{p_l}\}$ of $\{R_p\}$ such that each $R_{p_j}, \ R_{p_k}$  and  $R_{p_l}$  omits $u ~\mbox{and}~ v,$  $v ~\mbox{and}~ w,$ and $u ~\mbox{and}~ w,$ respectively. Without loss of generality, we assume $\{R_{p_j}\}$ to be a subsequence of $\{R_p\}$ which we we again denote by $R_p$ itself. Thus  each $R_p$ omits $u ~\mbox{and}~ v.$ Suppose $\{R_p\}$ is not normal at $z_0.$ Then by Zalcman's lemma we can find a sequence in $\{R_p\}$ which we again denote by $\{R_p\}$, a sequence $\{z_p\}$ of complex numbers with $z_p\to z_0$ and a sequence $\{\rho_p\}$ of positive real numbers with $\rho_p\to 0$ such that $S_p(\zeta)= R_p(z_p+\rho_p\zeta)$ converges locally uniformly with respect to the spherical metric to a non-constant meromorphic function $S$ on $\mathbb{C}.$ Since each $R_p$ omits $u ~\mbox{and}~ v,$ it follows that $S$ also omits $u ~\mbox{and}~ v$ and so $S$ takes some value $a$ infinitely many times. Let $\zeta_1,\zeta_2,\ldots ,\zeta_{m+1}$ be $m+1$ distinct zeros of $S-a.$ Since $S$  is non-constant, by Hurwitz's theorem there is a sequence $\zeta_{p,i}\to \zeta_i $ such that $R_p(z_p+\rho_p\zeta_{p,i})= a ~(i= 1,2,\ldots,m+1).$ Since $R_p$ is a rational function of degree at most $m,$  there is some $~i\neq j, 1\leq i,j\leq m+1$ such that $z_p+\rho_p\zeta_{p,i}= z_p+\rho_p\zeta_{p,j}$ which implies $\zeta_i=\zeta_j,$ a contradiction. Hence Theorem \ref{thmM} is valid for $\mathcal{R}_m.$\\
Further, as an immediate consequence of preceding discussion we find that if each $R\in \mathcal{R}_m.$ omits two given fixed values in $\mathbb{C}_{\infty},$ on a domain $D$, then $\mathcal{R}_m$ is normal in $D.$

\bibliographystyle{amsplain}

\end{document}